\documentclass[11pt]{article}
\usepackage{amsmath}
\usepackage{amssymb}
\usepackage{amsfonts}
\usepackage{color}
\allowdisplaybreaks

\newtheorem{Theorem}{\textbf{Theorem}}[section]
\newtheorem{Lemma}{\textbf{Lemma}}[section]
\newtheorem{Proposition}{\textbf{Proposition}}[section]
\newtheorem{Corollary}{\textbf{Corollary}}[section]
\newtheorem{Remark}{\textbf{Remark}}[section]
\newtheorem{Example}{\textbf{Example}}[section]
\newtheorem{Definition}{\textbf{Definition}}[section]

\newenvironment{theorem}{\begin{Theorem}$\!\!\!$}{\end{Theorem}}
\newenvironment{lemma}{\begin{Lemma}$\!\!\!$}{\end{Lemma}}
\newenvironment{proposition}{\begin{Proposition}$\!\!\!$}{\end{Proposition}}

\newenvironment{remark}{\begin{Remark}$\!\!\!$}{\end{Remark}}

\numberwithin{equation}{section}


\begin{document}
\title{$L^2$ asymptotic profiles of solutions to\\ 
linear damped wave equations}
\author{Hironori Michihisa${}^\ast$\\ 
{\small Department of Mathematics, Graduate School of Science, Hiroshima University} \\
{\small Higashi-Hiroshima 739-8526, Japan}
}

%
\date{}

\maketitle

\begin{abstract}
In this paper we obtain higher order asymptotic profilles of solutions to the Cauchy problem of the linear damped wave equation in $\textbf{R}^n$
\begin{equation*}
u_{tt}-\Delta u+u_t=0,
\qquad
u(0,x)=u_0(x),
\quad
u_t(0,x)=u_1(x), 
\end{equation*}
where $n\in\textbf{N}$ and $u_0$, $u_1\in L^2(\textbf{R}^n)$. 
Established hyperbolic part of asymptotic expansion seems to be new in the sense that the order of the expansion of the hyperbolic part depends on the spatial dimension.
\end{abstract}

\footnote[0]{\hspace{-2em} ${}^\ast$Corresponding author.}
\footnote[0]{\hspace{-2em} \textit{Email:} hi.michihisa@gmail.com}
\footnote[0]{\hspace{-2em} 2010 \textit{Mathematics Subject Classification}. 35B40, 35L05, 35L15}
\footnote[0]{\hspace{-2em} \textit{Keywords and Phrases}: Damped wave equation, Asymptotic expansion, Diffusion phenomena, Wave phenomena}

\section{Introduction}
In this paper we consider the Cauchy problem of the solution to the linear damped wave equation 
\begin{equation}
\label{eq:LDW}
\begin{cases}
u_{tt}-\Delta u+u_t=0, & t>0,\quad x\in\textbf{R}^n, \\
u(0,x)=u_0(x), & x\in\textbf{R}^n, \\
u_t(0,x)=u_1(x),& x\in\textbf{R}^n, \\
\end{cases}
\end{equation}
where $n\in\textbf{N}$ and $u_0$, $u_1\in L^2(\textbf{R}^n)$ are given initial data.  

Our purpose is to obtain higher order asymptotic expansions of the solution to problem \eqref{eq:LDW} in the $L^2$ framework for all spatial dimension $n$. 
We investigate both the wave and the diffusive structure. 
As will be described later in detail, we propose a quite simple idea to obtain such profiles by applying  the Taylor theorem to the appropriate functions, which are naturally included in the damped wave equation itself.

Now let us recall previous studies on the asymptotic behavior of solutions to damped wave equations. 

In \cite{M} Matsumura, who made a pioneering study on this subject, first obtained $L^m$-$L^2$ decay estimates of solutions. 
Since then, many mathematicians have analyzed the diffusion structure of hyperbolic type equations and systems (see e.g., \cite{HL}, \cite{HL2} and \cite{YM}).
His study in \cite{M} is also applied to the Cauchy problem of the nonlinear equation 
\begin{equation}
\label{eq:SLDW}
\begin{cases}
u_{tt}-\Delta u+u_t=|u|^p, & t>0,\quad x\in\textbf{R}^n, \\
u(0,x)=u_0(x), & x\in\textbf{R}^n, \\
u_t(0,x)=u_1(x),& x\in\textbf{R}^n,
\end{cases}
\end{equation}
where $p>1$. 
Concerning the problem \eqref{eq:SLDW}, 
Todorova-Yordanov \cite{ToYo} studied the so-called critical exponent problem, which implies that there exists a exponent $p^*>1$ of the nonlinear term $|u|^p$ such that if the power $p$ of the nonlinearity satisfies $p^*<p$, then the corresponding problem \eqref{eq:SLDW} has a unique small data global solution, while if $1<p\le p^*$, then the associated solution of \eqref{eq:SLDW} does not exist globally in time. 
In \cite{ToYo}, they derived $p^*=1+2/n$, which coincides with the Fujita exponent in the heat equation case \cite{F}. 
Independently of them \cite{ToYo}, Ikehata-Miyaoka-Nakatake \cite{IMN} succeeded in proving the global existence and optimal decay estimates of the total energy of the weak solutions to problem \eqref{eq:SLDW} with the power $p>p^*$ for $n=1$, $2$, and $p>2$ for $n=3$. 
It should be mentioned that the results in \cite{ToYo} fully depends on the compactness assumption on the support of initial data, 
while in \cite{IMN} they removed it. 
(See also \cite{K}, \cite{KNO}, \cite{IT}, \cite{Ni} and references therein.)
In this line, one should cite an extremely important result due to Nishihara \cite{Ni}, who derived the $L^p$-$L^q$ estimates for the difference between the solution $u(t,x)$ to problem \eqref{eq:LDW} and the sum of solutions $w(t,x)$, $\tilde{w}(t,x)$, and $v(t,x)$ corresponding to the equations 
\begin{equation*}
\label{eq:Wave}
\begin{cases}
w_{tt}-\Delta w=0, & t>0,\quad x\in\textbf{R}^n, \\
w(0,x)=u_0(x), & x\in\textbf{R}^n, \\
w_t(0,x)=u_1(x),& x\in\textbf{R}^n,
\end{cases}
\end{equation*}
\begin{equation*}
\label{eq:Wave2}
\begin{cases}
\tilde{w}_{tt}-\Delta \tilde{w}=0, & t>0,\quad x\in\textbf{R}^n, \\
\tilde{w}(0,x)=0, & x\in\textbf{R}^n, \\
\tilde{w}_t(0,x)=u_0(x),& x\in\textbf{R}^n,
\end{cases}
\end{equation*}
and heat equation
\begin{equation*}
\label{eq:Heat}
\begin{cases}
v_t-\Delta v=0, & t>0,\quad x\in\textbf{R}^n, \\
v(0,x)=u_0(x)+u_1(x), & x\in\textbf{R}^n, 
\end{cases}
\end{equation*}
respectively. 
In fact, Nishihara in \cite{Ni} derived 
\begin{equation}
\label{as}
u\sim
e^{-\frac{t}{2}}\left\{
w
+\left(
\frac{1}{2}+\frac{t}{8}
\right)\tilde{w}
\right\}
+v,
\qquad
t\to\infty,
\end{equation}
and this work has been done for $n=3$. 
Before \cite{Ni}, to the best of the author's knowledge, no attempt has been done to derive such a wave effect appeared in the first term of the right-hand side of  \eqref{as}. 
This wave effect has been considered as an error term for it decays much faster that diffusive part $v$. 
When Nishihara obtains and evaluates the wave effect, he suggested a method to divide the solution $u$ into two parts: 
one is a wave part multiplied by $e^{-\frac{t}{2}}$, 
and the other one is an integrated function whose integrand includes the modified Bessel function. 
Nowadays, this fashion is known as the Nishihara decomposition. 
Narazaki \cite{Na} obtained the $L^p$-$L^q$ estimates for general $n$, however its asymptotic profile increases as the spatial dimension $n$ does. 
This fact is already pointed out by Hosono-Ogawa \cite{HO}.  
In \cite{HO}, they also showed the $L^p$-$L^q$ estimates for the same function as that of \cite{Ni} in the case of $n=2$. 
Marcati-Nishihara \cite{MN} dealt with the case $n=1$, and derived that the solution $u$ to problem \eqref{eq:LDW} with $n=1$ behaves like 
\begin{equation}
\label{as2}
u\sim
e^{-\frac{t}{2}}w+v,
\qquad
t\to\infty.
\end{equation}
It should be mentioned that in \eqref{as2} the term $e^{-\frac{t}{2}}(1/2+t/8)\tilde{w}$ is dropped as compared with \eqref{as}.
One cannot find a clear explanation why there is such a difference of the profiles between \eqref{as} and  \eqref{as2}. 
In this paper we will study the reason concerning the difference, and the work will be done for all $n$ (see Proposition~\ref{prop:1} and \ref{prop:2} below). 
This is one of our novelties. 

Quite recently, Sakata-Wakasugi \cite{SW} extended the idea of \cite{CH} and \cite{Ni} to obtain the asymptotic profile of the solution to \eqref{eq:LDW} for general $n\in\textbf{N}$, and applied it to the study of the hot spots of solutions to \eqref{eq:LDW}.

Many asymptotic expansions of solutions to parabolic type equations are well-studied 
(see \cite{EZ}, \cite{HKN}, \cite{HKN2}, \cite{IIK}, \cite{IKK}, \cite{IKM} and references therein). 
From the viewpoint of the diffusion phenomena of the solution to \eqref{eq:LDW}, 
it is quite natural to consider the higher order expansions to the solution to problem \eqref{eq:LDW}. 
Takeda~\cite{T} recently succeeded in deriving the asymptotic profile of diffusive part of the solution to \eqref{eq:LDW}, and also obtained $L^p$-$L^q$ estimates together with weighted $L^1$ estimates. 
In his study the following asymptotic profile is obtained; 
for arbitrary $m\in\textbf{N}$,   
\begin{equation}
\begin{split}
\label{eq:takeda}
\hat{u}(t,\xi) 
& \sim 
\frac{1}{2}\sum_{j=0}^m \sum_{k=0}^{m-j} 
\alpha_{j,k} (-t)^j |\xi|^{2(2j+k)}e^{-t|\xi|^2}\widehat{u_0} \\
& \qquad\qquad 
+\sum_{j=0}^m \sum_{k=0}^{m-j} \sum_{\ell=0}^{m-j-k} 
\alpha_{j,k}\beta_\ell 
(-t)^j |\xi|^{2(2j+k+\ell)}e^{-t|\xi|^2}\widehat{u_1}.
\end{split}
\end{equation}
Here the coefficients $\alpha_{j,k}$ and $\beta_\ell$ are 
\begin{equation}
\label{coe}
\alpha_{j,k}
:=\frac{1}{j!k!}
\left.\frac{d^k}{dr^k}\phi_j(r)\right|_{r=0}, 
\qquad
\beta_\ell
:=\frac{1}{\ell!}
\left.\frac{d^\ell}{dr^\ell}\psi(r)\right|_{r=0},
\end{equation}
where 
\begin{equation*}
\label{pp}
\phi_j(r):=\left(\frac{1}{1/2+\sqrt{1/4-r}}\right)^{2j},
\qquad
\psi(r):=\frac{1}{2\sqrt{1/4-r}}.
\end{equation*}
The method to obtain \eqref{eq:takeda} is called the Takeda expansion. 
The number of summation symbols in \eqref{eq:takeda} represents how many times the Taylor theorem is applied. 
As a result, the resulting expansion \eqref{eq:takeda} is so complicated and redundant in the sense that some of coefficients are $0$, for example, $\alpha_{0,k}=0$ for $k\in\textbf{N}$. 

As will be seen later, this complexity can be avoided by discovering appropriate functions to which the Taylor theorem can be applied. 
The appropriate functions are naturally determined by the equation itself. 
The diffusive part of the asymptotic profile developed in this paper seems to be different from \eqref{eq:takeda}, but it will be shown that it actually coincides with \eqref{eq:takeda} (see Appendix). 
\\

In the rest of this section we explain our main idea to obtain precise wave and diffusive structures. 
As is often explained, the damped wave equation is a cross between the wave equation (when $u_t$ drops), and the heat equation (if $u_{tt}$ is removed).
So the solution to \eqref{eq:LDW} has wave-like property and behaves like the heat flow. 
From this observations, when one wants to see the diffusive structure, one considers the equation  
\begin{equation}
\label{eq:heat-v}
au_{tt}-\Delta u+u_t=0
\end{equation}
with small $a>0$. 
Then, in order to describe the solution of \eqref{eq:LDW} (=\eqref{eq:heat-v} with a=1), 
we first solve \eqref{eq:heat-v} in the Fourier space.  
In fact, the solution $\hat{u}(t,\xi;a)$ of \eqref{eq:heat-v} is expressed by 
\begin{equation}
\label{eq:SOL}
\hat{u}(t,\xi;a)
=e^{-\frac{t}{2a}}\cos\left(\frac{t\sqrt{a|\xi|^2-\frac{1}{4}}}{a}\right)\widehat{u_0} 
+e^{-\frac{t}{2a}}\frac{\sin\left(\frac{t\sqrt{a|\xi|^2-\frac{1}{4}}}{a}\right)}{\sqrt{a|\xi|^2-\frac{1}{4}}}\left(\frac{1}{2}\widehat{u_0}+a\widehat{u_1}\right).
\end{equation}
Although $a=0$ may be a singular limit, it works well in the diffusive part. 
For example, when one divides the first term in the right-hand side of \eqref{eq:SOL} in the low frequency parts as follows 
\begin{equation*}
e^{-\frac{t}{2a}}\cos\left(\frac{t\sqrt{a|\xi|^2-\frac{1}{4}}}{a}\right)
=\frac{1}{2}\exp\left(-\frac{2t|\xi|^2}{1+\sqrt{1-4a|\xi|^2}}\right)
-\frac{1}{2}e^{-t\left(\frac{1}{2a}+\sqrt{1-4a|\xi|^2}\right)}, 
\qquad
|\xi|\ll1
\end{equation*}
(see also the proof of Proposition~\ref{prop:3}),  
one finds the first term on the right-hand side can be defined at $a=0$, 
and the heat kernel $e^{-t|\xi|^2}$ automatically appears.
By this consideration, we naturally encounter nice functions $g$ and $h$, and similarly $f$ to be expanded asymptotically (see \eqref{eq:D-func} below). 

It will be considered whether this method works or not in such ``hybrid" type equations. 
The strongly damped wave equation 
\begin{equation*}
u_{tt}-\Delta u-\Delta u_t=0
\end{equation*}
will be studied in a forthcoming paper.
\\

The rest of this paper is organized as follows. 
In section~2 we describe the asymptotic profile of the solution to \eqref{eq:LDW} and state the main theorem. 
In section~3 we prepare several lemmas for the proof of Theorem~\ref{thm:1}. 
In section~4 we give the proof of Theorem~\ref{thm:1}, which directly follows from the results in section~3. 
As in Appendix, we will see that the asymptotic profile obtained in this paper includes the results in \cite{HO}, \cite{MN}, \cite{Ni}, and \cite{T}. 
\section{Main Result}
The solution $u=u(t,x)$ of \eqref{eq:LDW} is given by (see e.g., \cite{Ni})
\begin{equation*}
u(t)=K_0(t)u_0+K_1(t)\left(\frac{1}{2}u_0+u_1\right),
\end{equation*}
where 
\begin{equation*}
\label{eq:mul}
K_0(t)g
=\mathcal{F}^{-1}\left[e^{-\frac{t}{2}}\cos\left(t\sqrt{|\xi|^2-\frac{1}{4}}\right)\hat{g}\right]
=\mathcal{F}^{-1}\left[e^{-\frac{t}{2}}\cosh\left(t\sqrt{\frac{1}{4}-|\xi|^2}\right)\hat{g}\right],
\end{equation*}
\begin{equation*}
\label{eq:mul2}
K_1(t)g
=\mathcal{F}^{-1}\left[e^{-\frac{t}{2}}\frac{\sin\left(t\sqrt{|\xi|^2-\frac{1}{4}}\right)}{\sqrt{|\xi|^2-\frac{1}{4}}}\hat{g}\right]
=\mathcal{F}^{-1}\left[e^{-\frac{t}{2}}\frac{\sinh\left(t\sqrt{\frac{1}{4}-|\xi|^2}\right)}{\sqrt{\frac{1}{4}-|\xi|^2}}\hat{g}\right].
\end{equation*}
Here the function $\hat{g}$ represents the Fourier transform of $g$, 
\begin{equation*}
\hat{g}(\xi)=\int_{\textbf{R}^n}e^{-ix\cdot\xi}g(x)\,dx.
\end{equation*}
The inverse Fourier transform $\mathcal{F}^{-1}$ is similarly defined. 

We introduce cut-off functions. 
Let $\chi_L$, $\chi_H$ and $\chi_M$ be smooth functions such that 
\begin{equation*}
\chi_L(r) :=
\begin{cases}
0, & r\ge\frac{1}{3}, \\
1, & r\le\frac{1}{4}, 
\end{cases}
\qquad
\chi_H(r) :=
\begin{cases}
1, & r\ge2, \\
0, & r\le1, 
\end{cases}
\qquad 
\chi_M(r)  :=1-\chi_L(r)-\chi_H(r).
\end{equation*}
Set 
\begin{equation}
\label{eq:D-func}
\begin{split}
& f(r,c,t):=\cos\left(t\sqrt{r^2-c}\right), \\
& g(r,a,t):=\exp\left(-\frac{2tr^2}{1+\sqrt{1-4ar^2}}\right), \\
& h(r,a,t):=\frac{1}{\sqrt{1-4ar^2}}\exp\left(-\frac{2tr^2}{1+\sqrt{1-4ar^2}}\right).
\end{split}
\end{equation}
In this paper $\textbf{N}$ denotes the set of all natural numbers and we set $\textbf{N}_0:=\textbf{N}\cup\{0\}$. 

Let $b\in\textbf{N}$ and $\ell\in\textbf{N}$. 
Define 
\begin{equation*}
W^1_b(\xi,t):=
\sum_{k=0}^{b-1} \left(\frac{1}{4}\right)^k \frac{1}{k!} \frac{\partial^k f}{\partial c^k}\left(|\xi|,0,t\right),
\end{equation*}
\begin{equation*}
D^1_\ell(\xi,t):=
\frac{1}{2}\sum_{k=0}^{\ell-1} \frac{1}{k!} \frac{\partial^k g}{\partial a^k}(|\xi|,0,t),
\end{equation*}
and 
\begin{equation*}
W^2_b(\xi,t):=
2t^{-1}\sum_{k=0}^{b-2} \left(\frac{1}{4}\right)^k 
\frac{1}{k!} \frac{\partial^{k+1} f}{\partial c^{k+1}}\left(|\xi|,0,t\right),
\end{equation*}
i.e.,
\begin{equation*}
W^2_b(\xi,t):=
\begin{cases}
\displaystyle{2t^{-1}\sum_{k=0}^{b-2} \left(\frac{1}{4}\right)^k 
\frac{1}{k!} \frac{\partial^{k+1} f}{\partial c^{k+1}}\left(|\xi|,0,t\right)}, & b\ge2, \\
0, & b=1,
\end{cases}
\end{equation*}
\begin{equation*}
D^2_\ell(\xi,t):=
\sum_{k=0}^{\ell-1} \frac{1}{k!} \frac{\partial^k h}{\partial a^k}(|\xi|,0,t),
\end{equation*}
for $\xi\in\textbf{R}^n$ and $t>0$.
We also define 
\begin{equation}
\label{eq:m}
\begin{split}
m_{b,\ell,J}^1(\xi,t):=\chi_J(|\xi|)
\left[
e^{-\frac{t}{2}}\cos\left(t\sqrt{|\xi|^2-\frac{1}{4}}\right)
-e^{-\frac{t}{2}}W^1_b(\xi,t)
-D^1_\ell(\xi,t)
\right],
\end{split}
\end{equation}
\begin{equation}
\label{eq:m2}
\begin{split}
m_{b,\ell,J}^2(\xi,t):=\chi_J(|\xi|)
\left[
e^{-\frac{t}{2}}\frac{\sin\left(t\sqrt{|\xi|^2-\frac{1}{4}}\right)}{\sqrt{|\xi|^2-\frac{1}{4}}}
-e^{-\frac{t}{2}}W^2_b(\xi,t)
-D^2_\ell(\xi,t)
\right],
\end{split}
\end{equation}
for $J=L$, $M$, $H$. 

For $f\in L^2(\textbf{R}^n)$, we set 
\begin{equation*}
\|f\|_2
:=\left(
\int_{\textbf{R}^n} |f(x)|^2\, dx
\right)^{\frac{1}{2}}.
\end{equation*}

We are ready to state the main theorem. 
\begin{theorem}
\label{thm:1}
Let $n\in\textbf{N}$, $b\in\textbf{N}$ with $b>n/2$, and $\ell\in\textbf{N}$ 
and let $u$ be the solution of \eqref{eq:LDW}. 
Then there exists a constant $C>0$ such that 
\begin{equation*}
\begin{split}
& \left\|\hat{u}(t)
-\left(e^{-\frac{t}{2}}W_b^1(t)
+D_\ell^1(t)
\right)
\widehat{u_0}
-\left(e^{-\frac{t}{2}}W_b^2(t)
+D_\ell^2(t)
\right)
\left(\frac{1}{2}\widehat{u_0}+\widehat{u_1}\right)
\right\|_2 \\
& \le
\left(Ct^{2(b-1)}e^{-\frac{t}{2}}+Ct^b e^{-\frac{t}{2}}+Ct^{-\frac{n}{4}-\ell}\right)\|u_0\|_2 \\
& \qquad\qquad
+\left(Ct^{2(b-1)-1}e^{-\frac{t}{2}}+Ct^{b-1} e^{-\frac{t}{2}}+Ct^{-\frac{n}{4}-\ell}\right)\|u_1\|_2
\end{split}
\end{equation*}
for $t\ge1$. 
Here the constant $C>0$ does not depend on $t$ and the initial data $u_0$, $u_1\in L^2(\textbf{R}^n)$.
\end{theorem}
\begin{remark}
A discovery of the condition $b>n/2$ is our main contribution. 
The leading term of the hyperbolic part depends on the spatial dimension $n$. 
\end{remark}
\section{Preliminaries}
Now we recall Fa\`a di Bruno's formula which is a generalization of the chain rule to higher order derivatives.
Let $k\in\textbf{N}$, and let $F(x)$ and $G(x)$ be functions for which all necessary derivatives are defined, then 
\begin{equation*}
\label{eq:faa}
\frac{d^k}{dx^k}F(G(x))
=\sum_{p_1,\dots,p_k}^* 
\frac{k!}{\prod_{j=1}^k p_j! j!^{p_j}}
F^{(\sum_{j=1}^k p_j)}(G(x))
\prod_{j=1}^{k}\left(G^{(j)}(x)\right)^{p_j}.
\end{equation*}
Here and subsequently, the sum $\displaystyle{\sum_{p_1,\dots,p_k}^*}$ is taken over all 
$(p_1,\dots,p_k)\in{\textbf{N}_0}^k$ satisfying 
\begin{equation*}
\label{eq:dio}
\sum_{j=1}^k jp_j=k.
\end{equation*}

This formula is useful when we write down all the terms in higher order derivatives of a composite function in the sense that it  describes all the coefficients explicitly. 

\subsection{Estimates in High-frequency Region}

In this subsection, we obtain pointwise estimates for derivatives of $\cos\left(t\sqrt{r^2-c}\right)$ in the region $r\ge1$ and $0\le c\le1/4$. 
In this setting, we use the fact that 
\begin{equation*}
\sqrt{r^2-c}\ge \frac{\sqrt{3}}{2}r 
\qquad\mbox{if}\qquad
r\ge1
\quad\mbox{and}\quad
0\le c\le\frac{1}{4},
\end{equation*}
without mentioning anything afterwards. 

\begin{lemma}
\label{lem:cos}
Let $b\in\textbf{N}$. 
Then there exists a constant $C_b^1>0$ such that 
\begin{equation}
\label{ineq:M2}
\left|
\frac{\partial^b}{\partial c^b}\cos\left(t\sqrt{r^2-c}\right)
\right|
\le C_b^1 t^b r^{-b}
\end{equation}
for $t\ge1$, $r\ge1$ and $0\le c\le1/4$. 
In particular, there exists a constant $C_b^2>0$ such that 
\begin{equation}
\label{ineq:M3}
\left|\left.
\frac{\partial^b}{\partial c^b}\cos\left(t\sqrt{r^2-c}\right)
\right|_{c=0} \right|
\le C_b^2 t^b r^{-b}
\end{equation}
for $t\ge1$ and $r>0$.
\end{lemma}
\textbf{Proof.} 
It suffices to consider the case $b\ge2$. 
Inductively, we have
\begin{equation*}
\frac{\partial^j}{\partial c^j}
\left(\sqrt{r^2-c}\right)
=L_j\left(\sqrt{r^2-c}\right)^{-(2j-1)}.
\end{equation*}
for $j\in\textbf{N}$. 
Here and after, $L_j$ stands for the constant such that 
\begin{equation}
\label{eq:const}
L_j:=
\begin{cases}
\displaystyle{-\frac{(2j-3)!!}{2^j}}, & j\ge2, \\[10pt]
\displaystyle{-\frac{1}{2}}, & j=1.
\end{cases}
\end{equation}
Set $G(r,c):=\sqrt{r^2-c}$.  
Then, for each $j\in\textbf{N}$, there exists a constant $C_j>0$ such that 
\begin{equation}
\label{ineq:ST2}
\left|
\frac{\partial^j G}{\partial c^j}(r,c)
\right|
\le C_j r^{-(2j-1)} 
\qquad\mbox{for}\qquad
r\ge1
\quad\mbox{and}\quad
0\le c\le1/4. 
\end{equation}
We also see that 
\begin{equation}
\label{eq:ST3}
\frac{\partial^j G}{\partial c^j}(r,0)
=L_jr^{-(2j-1)}
\qquad\mbox{for}\qquad
r>0.
\end{equation}
Applying Fa\`a di Bruno's formula, we have 
\begin{equation*}
\begin{split}
& \frac{\partial^b}{\partial c^b}\cos\biggr(tG(r,c)\biggr)
=\sum_{p_1,\dots,p_b}^* \frac{b!}{\prod_{j=1}^b p_j! j!^{p_j}} 
\cos^{\left(\sum_{j=1}^b p_j\right)} \biggr(tG(r,c)\biggr)
\prod_{j=1}^b \left(t\frac{\partial^j G}{\partial c^j}(r,c)\right)^{p_j} \\
& =t^b \left(\frac{\partial G}{\partial c}(r,c)\right)^b \cos^{(b)}\biggr(tG(r,c)\biggr) \\
& \qquad\qquad
+\sum_{p_1,\dots,p_b} \frac{b!}{\prod_{j=1}^b p_j! j!^{p_j}} 
t^{\sum_{j=1}^b p_j}
\cos^{\left(\sum_{j=1}^b p_j\right)} \biggr(tG(r,c)\biggr)
\prod_{j=1}^b \left(\frac{\partial^j G}{\partial c^j}(r,c)\right)^{p_j},
\end{split}
\end{equation*}
where the sum $\sum_{p_1,\dots,p_b}$ is taken over all $(p_1,\dots,p_b)\in{\textbf{N}_0}^b$ satisfying 
\begin{equation*}
\sum_{j=1}^b jp_j=b
\quad
\mbox{with}
\quad
p_1\le b-1. 
\end{equation*}
It follows from \eqref{ineq:ST2} that 
\begin{equation*}
\begin{split}
\left|
\frac{\partial^b}{\partial c^b}\cos\left(t\sqrt{r^2-c}\right)
\right|
& \le Ct^b r^{-b}
+C\sum_{p_1,\dots,p_b}
t^{\sum_{j=1}^b}
r^{-\sum_{j=1}^b (2j-1)p_j} \\
& = Ct^b r^{-b}
+C\sum_{p_1,\dots,p_b}
t^{\sum_{j=1}^b}
r^{-2b+\sum_{j=1}^b p_j} \\
& \le Ct^b r^{-b}
+C\sum_{p_1,\dots,p_b}
t^{b-1}
r^{-2b+(b-1)} \\
& =Ct^b r^{-b}
+Ct^{b-1}r^{-b-1} 
\le Ct^b r^{-b}
\end{split}
\end{equation*}
for $t\ge1$, $r\ge1$ and $0\le c\le1/4$.
Hence, we obtain \eqref{ineq:M2}. 
Inequality \eqref{ineq:M3} is proved in the same way as above 
with the aid of \eqref{eq:ST3}. 
Therefore the proof is complete. 
$\Box$
%

\subsection{Estimates in Low-frequency Region}
Here we often use the relation 
\begin{equation*}
\sqrt{1-4ar^2}\ge \frac{\sqrt{5}}{3}
\qquad\mbox{if}\qquad
0\le r\le\frac{1}{3}
\quad\mbox{and}\quad
0\le a\le1.
\end{equation*}

To begin with, we check the asymptotic profile of wave part (see Appendix) is well-defined as an $L^2$ function in the low-frequency region. 
To do this, our first goal is to prove Lemma~\ref{lem:SING}. 
For this purpose, we prepare several lemmas. 
\begin{lemma}
\label{lem:kinou}
Define 
\begin{equation*}
F_k(r,c,t):=
\frac{\partial^k}{\partial c^k} 
\cos\left(t\sqrt{r^2-c}\right)
\end{equation*}
for $k\in\textbf{N}_0$.
Then 
\begin{equation}
\label{eq:IND}
F_k(r,c,t)=
\frac{1}{4(r^2-c)}
\biggr(-t^2 F_{k-2}(r,c,t)+2(2k-3)F_{k-1}(r,c,t)\biggr)
\end{equation}
for $k\ge2$.
\end{lemma}
\textbf{Proof.} 
The proof is done by induction on $k$. 
Throughout the proof, we write $F_k=F_k(r,c,t)$ for simplicity. 
First we calculate 
\begin{equation*}
F_1=\frac{t}{2\sqrt{r^2-c}}\sin\left(t\sqrt{r^2-c}\right),
\end{equation*}
\begin{equation*}
\begin{split}
F_2
& =-\frac{t^2}{4(r^2-c)}\cos\left(t\sqrt{r^2-c}\right)
+\frac{t}{4\left(\sqrt{r^2-c}\right)^3}\sin\left(t\sqrt{r^2-c}\right) \\
& =-\frac{t^2}{4(r^2-c)}F_0
+\frac{1}{2(r^2-c)}F_1 
=\frac{1}{4(r^2-c)}(-t^2 F_0+2F_1),
\end{split}
\end{equation*}
which shows that the statement is true for $k=2$.

Assuming \eqref{eq:IND} to hold for $k\ge2$, 
we will prove it for $k+1$.
It follows that 
\begin{equation*}
\begin{split}
F_{k+1}
& =\frac{1}{4(r^2-c)^2}
\biggr(-t^2 F_{k-2}+2(2k-3)F_{k-1}\biggr)
+\frac{1}{4(r^2-c)}
\biggr(-t^2 F_{k-1}+2(2k-3)F_k\biggr) \\
& =\frac{1}{r^2-c}F_k 
+\frac{1}{4(r^2-c)}
\biggr(-t^2 F_{k-1}+2(2k-3)F_k\biggr) \\
& =\frac{1}{4(r^2-c)}
\biggr(4F_k-t^2 F_{k-1}+2(2k-3)F_k\biggr) \\
& =\frac{1}{4(r^2-c)}
\biggr(-t^2 F_{k-1}+2(2k-1)F_k\biggr).
\end{split}
\end{equation*}
Thus \eqref{eq:IND} is also true for $k+1$. 
Therefore the proof is complete. $\Box$

From \eqref{eq:IND} of Lemma~\ref{lem:kinou} we see that  
\begin{equation*}
\begin{split}
F_k(r,0,t)
& =\frac{1}{4}\frac{-t^2 F_{k-2}(r,0,t)+2(2k-3)F_{k-1}(r,0,t)}{r^2} \\
& =\frac{1}{4}\frac{\biggr(-t^2 F_{k-2}(r,0,t)+2(2k-3)F_{k-1}(r,0,t)\biggr)r^{2k-3}}{r^{2k-1}}.
\end{split}
\end{equation*}
Set 
\begin{equation*}
I_k(r,t):=
\begin{cases}
\displaystyle{\frac{1}{2}t\sin(tr)}, & k=1,\\
\\
\displaystyle{\frac{1}{4}\biggr(-t^2 F_{k-2}(r,0,t)+2(2k-3)F_{k-1}(r,0,t)\biggr)r^{2k-3}}, & k\ge2.
\end{cases}
\end{equation*}
Then, one has 
\begin{equation*}
F_k(r,0,t)=\frac{I_k(r,t)}{r^{2k-1}}
\end{equation*}
for $k\in\textbf{N}$.

\begin{lemma}
\label{lem:I1}
For all $k\in\textbf{N}$, it holds that  
\begin{equation*}
I_k(0,t)=0
\end{equation*}
for $t>0$.
\end{lemma}
\textbf{Proof.} 
Fa\`a di Bruno's formula and \eqref{eq:const} yield 
\begin{equation*}
\begin{split}
F_k(r,0,t)
& =\sum_{p_1,\dots,p_k}^* 
\frac{k!}{\prod_{j=1}^k p_j! j!^{p_j}}
\cos^{(\sum_{j=1}^k p_j)}(tr)
\prod_{j=1}^{k}\left(tL_j r^{-2j+1} \right)^{p_j} \\
& =\sum_{p_1,\dots,p_k}^* 
\frac{k!}{\prod_{j=1}^k p_j! j!^{p_j}}
\left(\prod_{j=1}^k L_j^{p_j}\right)
t^{\sum_{j=1}^k p_j}
r^{-2k+\sum_{j=1}^k p_j}
\cos^{(\sum_{j=1}^k p_j)}(tr) \\
& =-L_k tr^{-2k+1} \sin(tr) \\
& \qquad \qquad 
+\sum_{p_1,\dots,p_k} 
\frac{k!}{\prod_{j=1}^k p_j! j!^{p_j}}
\left(\prod_{j=1}^k L_j^{p_j}\right)
t^{\sum_{j=1}^k p_j}
r^{-2k+\sum_{j=1}^k p_j}
\cos^{(\sum_{j=1}^k p_j)}(tr),
\end{split}
\end{equation*}
where the sum $\displaystyle{\sum_{p_1,\dots,p_k}}$ is taken over all $(p_1,\dots,p_k)\in{\textbf{N}_0}^k$ satisfying 
\begin{equation*}
\sum_{j=1}^k jp_j=k
\quad
\mbox{with}
\quad
p_k=0.
\end{equation*}
This constraint shows 
\begin{equation*}
\sum_{j=1}^k p_j
=\sum_{j=1}^{k-1} p_j
\ge 2.
\end{equation*}
Substituting $r=0$ into the following equality
\begin{equation*}
\begin{split}
I_k(r,t)
& =r^{2k-1}F_k(r,0,t) \\
& =-L_k t\sin(tr) 
+\sum_{p_1,\dots,p_k} 
\frac{k!}{\prod_{j=1}^k p_j! j!^{p_j}}
\left(\prod_{j=1}^k L_j^{p_j}\right)
t^{\sum_{j=1}^k p_j}
r^{-1+\sum_{j=1}^k p_j}
\cos^{(\sum_{j=1}^k p_j)}(tr),
\end{split}
\end{equation*}
we finish the proof of lemma. 
$\Box$
\begin{lemma}
\label{lem:I2}
Let $2\le k\in\textbf{N}$. 
Then 
\begin{equation}
\label{eq:Dr}
\frac{1}{r}\frac{\partial}{\partial r}I_k(r,t)
=\frac{t^2}{2} I_{k-1}(r,t).
\end{equation}
\end{lemma}
\textbf{Proof.} 
The proof is done by induction on $k$. 
First we see that 
\begin{equation*}
\begin{split}
I_2(r,t)
& =-\frac{t^2}{4} r\cos(tr)+\frac{t}{4}\sin(tr),
\end{split}
\end{equation*}
\begin{equation*}
\begin{split}
I_3(r,t)
& =-\frac{t^3}{8}r^2\sin(tr)
-\frac{3t^2}{8}r\cos(tr)
+\frac{3t}{8}\sin(tr)
\end{split}
\end{equation*}
and so we have 
\begin{equation*}
\begin{split}
I'_2(r,t)
 =\frac{t^3}{4} r\sin(tr),
\qquad 
I'_3(r,t)
& =-\frac{t^4}{8} r^2\cos(tr)
+\frac{t^3}{8} r\sin(tr).
\end{split}
\end{equation*}
Here we use symbol $'$ as $\partial/\partial r$ for simplicity.
Thus it follows that 
\begin{equation*}
I'_2(t,r)=\frac{t^2 r}{2}I_1(t,r),
\qquad 
I'_3(t,r)=\frac{t^2 r}{2}I_2(t,r).
\end{equation*}
Hence \eqref{eq:Dr} is true for $k=2$, $3$. 

Next assuming \eqref{eq:Dr} holds for $k$ and $k-1$ with $k\ge3$, 
we will prove it for $k+1$.
By definition of $I_k$ we have
\begin{equation*}
I_{k+1}(t,r)
=\frac{2k-1}{2}I_k(t,r)
-\frac{t^2}{4}r^2 I_{k-1}(t,r).
\end{equation*}
Thus we obtain 
\begin{equation*}
\begin{split}
\frac{1}{r} I'_{k+1}(t,r)
& =\frac{2k-1}{2}\frac{1}{r}I'_k(t,r)
-\frac{t^2}{4}r^2 \frac{1}{r}I'_{k-1}(t,r)
-\frac{t^2}{2} I_{k-1}(t,r) \\
& =\frac{2k-1}{2}\frac{t^2}{2}I_{k-1}(t,r)
-\frac{t^2}{4}r^2 \frac{t^2}{2}I_{k-2}(t,r)
-\frac{t^2}{2} I_{k-1}(t,r) \\
& =\frac{t^2}{2}\left(
\frac{2k-3}{2}I_{k-1}(t,r) 
-\frac{t^2}{4}r^2 I_{k-2}(t,r)
\right) \\
& =\frac{t^2}{2}I_k(t,r),
\end{split}
\end{equation*}
and so \eqref{eq:Dr} is true for $k+1$. 
Therefore the proof is now complete. 
$\Box$

\begin{lemma}
\label{lem:SING}
Let $k\in\textbf{N}$. 
Then 
\begin{equation}
\label{eq:SING}
\lim_{r\downarrow0}\left(
\left.\frac{\partial^k}{\partial c^k}\cos\left(t\sqrt{r^2-c}\right)\right|_{c=0}\right)
=\frac{t^{2k}}{2^k (2k-1)!!}
\end{equation}
for $t>0$.
\end{lemma}
\textbf{Proof.} 
By Lemmas~\ref{lem:I1} and \ref{lem:I2} and \eqref{eq:Dr}, applying the l'Hospital rule $k$-times, it follows that 
\begin{equation}
\begin{split}
\label{eq:SING}
\lim_{r\downarrow0}\left(
\left.\frac{\partial^k}{\partial c^k}\cos\left(t\sqrt{r^2-c}\right)\right|_{c=0}\right)
& =\lim_{r\downarrow0} F_k(r,0,t) 
=\lim_{r\downarrow0} \frac{I_k(r,t)}{r^{2k-1}} \\
& =\frac{1}{2k-1}\lim_{r\downarrow0} \frac{I'_k(r,t)}{r^{2k-2}} \\
& =\frac{1}{2k-1}\frac{t^2}{2}\lim_{r\downarrow0} \frac{I_{k-1}(r,t)}{r^{2k-3}} \\
& =\cdots \\
& =\frac{1}{(2k-1)\cdots5}\left(\frac{t^2}{2}\right)^{k-2}
\lim_{r\downarrow0} \frac{I_2(r,t)}{r^3} \\
& =\frac{1}{(2k-1)\cdots5\cdot3}\left(\frac{t^2}{2}\right)^{k-1}
\lim_{r\downarrow0} \frac{I_1(r,t)}{r} \\
& =\frac{t^{2k}}{2^k (2k-1)!!}
\end{split}
\end{equation}
for $t>0$. 
$\Box$
\\

The latter part of this subsection is devoted to obtain the estimates for diffusive part. 
\begin{lemma}
\label{lem:arx}
Let $k\in\textbf{N}_0$ and set $G(r,a)=\sqrt{1-4ar^2}$. 
Then there exists a constant $C_k^1>0$ such that 
\begin{equation}
\label{ineq:h0}
\left|
\frac{\partial^k}{\partial a^k}\biggr(1+G(r,a)\biggr)^{-1}
\right|
\le C_k^1 r^{2k} 
\end{equation}
for $0\le r\le1/3$ and $0\le a\le1$. 
Furthermore, there exists a constant $C_k^2\not=0$ such that 
\begin{equation}
\label{eq:h0}
\left.
\frac{\partial^k}{\partial a^k}\biggr(1+G(r,a)\biggr)^{-1}
\right|_{a=0}
=C_k^2 r^{2k}
\end{equation}
for $r\ge0$. 
Here constatns $C_k^1$ and $C_k^2$ are independent of $r$ and $a$. 
\end{lemma}
\textbf{Proof.} 
It suffices to show the lemma for $k\in\textbf{N}$.
Fa\`a di Bruno's formula gives 
\begin{equation*}
\begin{split}
& \frac{\partial^k}{\partial a^k}\biggr(1+G(r,a)\biggr)^{-1}
=\sum_{p_1,\dots,p_k}^* 
\frac{k!}{\prod_{j=1}^k p_j! j!^{p_j}}
\left.\frac{d^{\sum_{j=1}^k}}{ds^{\sum_{j=1}^k}}s^{-1}\right|_{s=1+G(r,a)}
\prod_{j=1}^k 
\left(\frac{\partial^j G}{\partial a^j}(r,a)\right)^{p_j}.
\end{split}
\end{equation*}
We see that 
\begin{equation}
\label{eq:2con}
\frac{\partial^j G}{\partial a^j}(r,a)
=4^j L_j r^{2j} G(r,a)^{-(2j-1)},
\end{equation}
where $L_j$ is already defined in \eqref{eq:const}.
Note that $G(r,a)^{-1}$ is bounded for $0\le r\le1/3$ and $0\le a\le1$, and that $G(r,0)=1$.
Thus it follows that 
\begin{equation*}
\begin{split}
\left|
\frac{\partial^k}{\partial a^k}\biggr(1+G(r,a)\biggr)^{-1}
\right|
\le C\sum_{p_1,\dots,p_k}^* 
\prod_{j=1}^k r^{2jp_j} 
\le Cr^{2\sum_{j=1}^k jp_j} 
=C r^{2k}
\end{split} 
\end{equation*}
for $0\le r\le1/3$ and $0\le a\le1$. 
We also see that 
\begin{equation*}
\begin{split}
\left.
\frac{\partial^k}{\partial a^k}\biggr(1+G(r,a)\biggr)^{-1}
\right|_{a=0} 
& =\sum_{p_1,\dots,p_k}^* 
\frac{k!}{\prod_{j=1}^k p_j! j!^{p_j}}
\left.\frac{d^{\sum_{j=1}^k}}{ds^{\sum_{j=1}^k}}s^{-1}\right|_{s=2}
\prod_{j=1}^k 
C_j^{p_j} r^{2jp_j} \\
& =\left(
\sum_{p_1,\dots,p_k}^* 
\frac{k!}{\prod_{j=1}^k p_j! j!^{p_j}}
\left.\frac{d^{\sum_{j=1}^k}}{ds^{\sum_{j=1}^k}}s^{-1}\right|_{s=2}
\prod_{j=1}^k 
C_j^{p_j}
\right) r^{2k}.
\end{split}
\end{equation*}
for $r\ge0$.
$\Box$
\begin{lemma}
\label{lem:H0}
Let $k\in\textbf{N}$. Then there exists a constant $C>0$ such that 
\begin{equation}
\label{ineq:H0}
\left|
\frac{\partial^k}{\partial a^k}\exp\left(-\frac{2tr^2}{1+\sqrt{1-4ar^2}}\right)
\right|
\le Cr^{2k}e^{-tr^2}\sum_{j=1}^k (tr^2)^j 
\end{equation}
for $0\le r\le1/3$, $0\le a\le 1$ and $t>0$.
Furthermore, there exist constants $C_j>0$ such that 
\begin{equation}
\label{eq:H0}
\left.
\frac{\partial^k}{\partial a^k}\exp\left(-\frac{2tr^2}{1+\sqrt{1-4ar^2}}\right)
\right|_{a=0}
=r^{2k}e^{-tr^2}\sum_{j=1}^k C_j (tr^2)^j 
\end{equation}
for $r\ge0$ and $t>0$.
Here constants $C$ and $C_j$ $(j=1,\dots,k)$ are independent of $r$, $a$ and $t$. 
\end{lemma}
\textbf{Proof.} 
Set $G(r,a):=\sqrt{1-4ar^2}$. 
Applying Fa\`a di Bruno's formula, we have 
\begin{equation*}
\begin{split}
& \frac{\partial^k}{\partial a^k}
\exp\left(-\frac{2tr^2}{1+G(r,a)}\right) \\
& =\sum_{p_1,\dots,p_k}^* 
\frac{k!}{\prod_{j=1}^k p_j! j!^{p_j}}
\exp\left(-\frac{2tr^2}{1+G(r,a)}\right)
\prod_{j=1}^k 
\left(\frac{\partial^j}{\partial a^j}\left(-\frac{2tr^2}{1+G(r,a)}\right)\right)^{p_j} \\
& =\exp\left(-\frac{2tr^2}{1+G(r,a)}\right)
\sum_{p_1,\dots,p_k}^* 
\frac{k!}{\prod_{j=1}^k p_j! j!^{p_j}}
(-2tr^2)^{\sum_{j=1}^k p_j} 
\prod_{j=1}^k
\left(\frac{\partial^j}{\partial a^j}
\biggr(1+G(r,a)\biggr)^{-1}\right)^{p_j}.
\end{split}
\end{equation*}
It follows from \eqref{ineq:h0} of Lemma~\ref{lem:arx} that 
\begin{equation*}
\begin{split}
\left|
\frac{\partial^k}{\partial a^k}
\exp\left(-\frac{2tr^2}{1+G(r,a)}\right)
\right|
& \le Ce^{-tr^2}\sum_{p_1,\dots,p_k}^* 
(tr^2)^{\sum_{j=1}^k p_j}
\prod_{j=1}^k 
r^{2jp_j} \\
& =Ce^{-tr^2}r^{2k}
\sum_{p_1,\dots,p_k}^* 
(tr^2)^{\sum_{j=1}^k p_j} \\
& \le Ce^{-tr^2}r^{2k}
\sum_{j=1}^k (tr^2)^j
\end{split}
\end{equation*}
for $0\le r\le1/3$, $0\le a\le 1$ and $t>0$.
Thus we obtain \eqref{ineq:H0}.

We see from \eqref{eq:h0} of Lemma~\ref{lem:arx} that  
\begin{equation*}
\begin{split}
& \left.
\frac{\partial^k}{\partial a^k}\exp\left(-\frac{2tr^2}{1+\sqrt{1-4ar^2}}\right)
\right|_{a=0} \\
& =e^{-tr^2}
\sum_{p_1,\dots,p_k}^* 
\frac{k!}{\prod_{j=1}^k p_j! j!^{p_j}}
(-2tr^2)^{\sum_{j=1}^k p_j} 
\prod_{j=1}^k
\left(C_j^2 r^{2j}\right)^{p_j} \\
& =e^{-tr^2}r^{2k}
\sum_{p_1,\dots,p_k}^* 
\frac{k!}{\prod_{j=1}^k p_j! j!^{p_j}}
(-2tr^2)^{\sum_{j=1}^k p_j} 
\prod_{j=1}^k
\left(C_j^2\right)^{p_j} \\
& =e^{-tr^2}r^{2k}
\sum_{j=1}^k C_j (tr^2)^j 
\end{split}
\end{equation*}
for $r\ge0$ and $t>0$. 
At last, we have replaced the constants with different ones. 
Hence, \eqref{eq:H0} hols and the proof is complete.
$\Box$
\begin{lemma}
\label{lem:H1}
Let $k\in\textbf{N}$. Then there exists a constant $C>0$ such that 
\begin{equation}
\label{ineq:H1}
\left|
\frac{\partial^k}{\partial a^k}\left(
\frac{1}{\sqrt{1-4ar^2}}
\exp\left(-\frac{2tr^2}{1+\sqrt{1-4ar^2}}\right)
\right)\right|
\le Cr^{2k}e^{-tr^2}\sum_{j=0}^k (tr^2)^j 
\end{equation}
for $0\le r\le1/3$, $0\le a\le 1$ and $t>0$. 
Furthermore, there exist constants $C_j$ such that 
\begin{equation}
\label{eq:H1}
\left.
\frac{\partial^k}{\partial a^k}\left(
\frac{1}{\sqrt{1-4ar^2}}
\exp\left(-\frac{2tr^2}{1+\sqrt{1-4ar^2}}\right)
\right)\right|_{a=0}
= Cr^{2k}e^{-tr^2}\sum_{j=0}^k (tr^2)^j 
\end{equation}
for $0\le r\le1/3$, $0\le a\le 1$ and $t>0$. 
Here constants $C$ and $C_j$ $(j=1,\dots,k)$ are independent of $r$, $a$ and $t$. 
\end{lemma}
\textbf{Proof.} 
Again we set $G(r,a)=\sqrt{1-4ar^2}$. 
By the Leibniz rule we see that 
\begin{equation*}
\begin{split}
& \frac{\partial^k}{\partial a^k} 
\left(
G(r,a)^{-1} \exp\left(-\frac{2tr^2}{1+G(r,a)}\right)
\right) \\
& \qquad\qquad 
=\sum_{j=0}^k \frac{k!}{(k-j)!j!}
\left(\frac{\partial^{k-j}}{\partial a^{k-j}}G(r,a)^{-1}\right)
\left(\frac{\partial^j}{\partial a^j}\exp\left(-\frac{2tr^2}{1+G(r,a)}\right)\right).
\end{split}
\end{equation*}
On the other hand, for each $i\in\textbf{N}$, one has 
\begin{equation}
\label{eq:cco}
\frac{\partial^i}{\partial a^i}G(r,a)^{-1}
=2^{i}(2i-1)!! r^{2i} G(r,a)^{-(2i+1)}.
\end{equation}
This implies that 
\begin{equation}
\label{ineq:G-inv}
\left|
\frac{\partial^i}{\partial a^i}G(r,a)^{-1}
\right|
\le Cr^{2i}
\end{equation}
and that 
\begin{equation}
\label{eq:G-inv}
\left.
\frac{\partial^i}{\partial a^i}G(r,a)^{-1}
\right|_{a=0}
=2^{i}(2i-1)!! r^{2i}.
\end{equation}
Therefore, it follows form \eqref{ineq:H0} and \eqref{ineq:G-inv} that
\begin{equation*}
\begin{split}
& \left|
\frac{\partial^k}{\partial a^k}\left(
\frac{1}{\sqrt{1-4ar^2}}
\exp\left(-\frac{2tr^2}{1+\sqrt{1-4ar^2}}\right)
\right)\right| \\
& \le 
\left|
\frac{\partial^k}{\partial a^k}
G(r,a)^{-1}
\right| e^{-tr^2}
+C\sum_{j=1}^k
\left|
\frac{\partial^{k-j}}{\partial a^{k-j}}G(r,a)^{-1}
\right| 
\left|
\frac{\partial^j}{\partial a^j}
\exp\left(-\frac{2tr^2}{1+G(r,a)}\right)
\right| \\
& \le 
Cr^{2k} e^{-tr^2}
+C\sum_{j=1}^k
r^{2(k-j)}
r^{2j}e^{-tr^2} \sum_{i=1}^j (tr^2)^i \\
& = 
Cr^{2k} e^{-tr^2}
+Cr^{2k} e^{-tr^2} 
\sum_{j=1}^k 
\sum_{i=1}^j (tr^2)^i \\
& \le  
Cr^{2k} e^{-tr^2}
+Cr^{2k} e^{-tr^2} 
\sum_{j=1}^k (tr^2)^j \\
& = 
Cr^{2k} e^{-tr^2} 
\sum_{j=0}^k (tr^2)^j
\end{split}
\end{equation*}
for $0\le r\le1/3$, $0\le a\le 1$ and $t>0$. 
Thus we obtain \eqref{ineq:H1}. 

Next, we see from \eqref{eq:H0} and \eqref{eq:G-inv} that 
\begin{align*}
& \left.
\frac{\partial^k}{\partial a^k}\left(
\frac{1}{\sqrt{1-4ar^2}}
\exp\left(-\frac{2tr^2}{1+\sqrt{1-4ar^2}}\right)
\right)\right|_{a=0} \\
& = \left.
\frac{\partial^k }{\partial a^k}G(r,a)^{-1}
\right|_{a=0} e^{-tr^2}
+\sum_{j=1}^k \frac{k!}{(k-j)!j!}
\left.
\frac{\partial^{k-j}}{\partial a^{k-j}}G(r,a)^{-1}
\right|_{a=0} 
\left.
\frac{\partial^j}{\partial a^j}
\exp\left(-\frac{2tr^2}{1+G(r,a)}\right)
\right|_{a=0} \\
& = D_k r^{2k} e^{-tr^2}
+\sum_{j=1}^k \frac{k!}{(k-j)!j!}
D_{k-j}r^{2(k-j)}
r^{2j}e^{-tr^2}\sum_{i=1}^j C_i (tr^2)^i \\
& = D_k r^{2k} e^{-tr^2}
+\sum_{j=1}^k \frac{k!}{(k-j)!j!}
D_{k-j}r^{2(k-j)}
r^{2j}e^{-tr^2}\sum_{i=1}^j C_i (tr^2)^i \\
& =D_k r^{2k} e^{-tr^2}
+\sum_{j=1}^k \frac{k!}{(k-j)!j!}
D_{k-j}r^{2k}e^{-tr^2}\sum_{i=1}^j C_i (tr^2)^i \\
& =D_k r^{2k} e^{-tr^2}
+r^{2k}e^{-tr^2}\sum_{j=1}^k E_j (tr^2)^j \\
& =r^{2k}e^{-tr^2}\sum_{j=0}^k E_j (tr^2)^j.
\end{align*}
Here we have replaced some constants from $D_j$ to $E_j$. 
Therefore the proof is complete. 
$\Box$
\\

The following lemma is found in, e.g., \cite{T} and \cite{TY}.
\begin{lemma}
\label{lem:N2}
Let $n\in\textbf{N}$, $1\le q\le \infty$ and $k\ge0$. Then there exists a constant $C>0$ such that 
\begin{equation*}
\||x|^k e^{-t|x|^2}\|_{L^q(\{x\in\textbf{R}^n:|x|\le1\})}\le C(1+t)^{-\frac{n}{2q}-\frac{k}{2}}.
\end{equation*}
for $t\ge0$.
\end{lemma}
\textbf{Proof.} 
In this paper, we use this lemma with $q=2$,  
so we check it only for $q=2$. 
It follows that 
\begin{equation*}
\begin{split}
\||x|^k e^{-t|x|^2}\|_{L^2(|x|\le1)}^2
& =\int_{|x|\le1}|x|^{2k} e^{-2t|x|^2}\,dx 
 \le C\int_0^1 r^{2k+n-1}e^{-2tr^2}\,dr \\
& \le C\int_0^1 r^{2k+n-1}e^{-(1+t)r^2}\,dr 
 \le C(1+t)^{-\frac{2k+n}{2}}\int_0^{\sqrt{t}}R^{2k+n-1}e^{-R^2}\,dR \\
& \le C(1+t)^{-\frac{n}{2}-k}
\end{split}
\end{equation*}
for $t\ge0$.   
The proof is now complete. 
$\Box$
%

\section{Poof of Theorem}
Let $u$ be the solution of \eqref{eq:LDW}. 
Then we see that 
\begin{equation*}
\begin{split}
& \left\|
\hat{u}(t)
-\left(e^{-\frac{t}{2}}W_b^1(\xi,t)
+D_\ell^1(\xi,t)
\right)
\widehat{u_0}
-\left(e^{-\frac{t}{2}}W_b^2(\xi,t)
+D_\ell^2(\xi,t)
\right)
\left(\frac{1}{2}\widehat{u_0}+\widehat{u_1}\right)
\right\|_2 \\
& \le \sum_{J=L,M,H}
\biggr\{
\|m^1_{b,\ell,J}\|_2 \|u_0\|_2
+\|m^2_{b,\ell,J}\|_2 \left(\frac{1}{2}\|u_0\|_2+\|u_1\|_2\right)
\biggr\},
\end{split}
\end{equation*}
where $m_{b,\ell,J}^i$ ($i=1$, $2$; $J=L$, $M$, $H$) are defined in \eqref{eq:m} and \eqref{eq:m2}.
Proof of Theorem~\ref{thm:1} can be divided into the following six propositions 
by treating each $m_{b,\ell,J}^i$ ($i=1$, $2$), which will be estimated in the low, middle, and high frequency region with $j=L$, $M$, and $H$, respectively. 

First we deal with the case $m_{b,\ell,H}^i$ ($i=1$, $2$). 
It is essential to estimate the norm of wave part $W_b^i$ ($i=1$, $2$) on the support of $\chi_H$, i.e., in the high frequency region.
\begin{proposition}
\label{prop:1}
Let $n\in\textbf{N}$, $b\in\textbf{N}$ with $b>n/2$, and $\ell\in\textbf{N}$. 
Then there exists a constant $C>0$ such that 
\begin{equation*}
\left\|m^1_{b,\ell,H}(t)\right\|_2
\le Ct^b e^{-\frac{t}{2}}
\end{equation*}
for $t\ge1$.
\end{proposition}
\textbf{Proof.} 
The Taylor theorem gives  
\begin{equation*}
\begin{split}
\cos\left(t\sqrt{|\xi|^2-\frac{1}{4}}\right)-W^1_b(\xi,t)
& =f\left(|\xi|,\frac{1}{4}\right)
-\sum_{k=0}^{b-1} \frac{1}{k!}\left(\frac{1}{4}\right)^k \frac{\partial^k f}{\partial c^k}\left(|\xi|,0\right) \\
& =\frac{1}{b!}\left(\frac{1}{4}\right)^b \frac{\partial^b f}{\partial c^b}\left(|\xi|,\frac{1}{4}\tau \right)
\end{split}
\end{equation*}
for some $0\le\tau\le1$. 
It follows from \eqref{ineq:M2} that 
\begin{equation*}
\left|
\cos\left(t\sqrt{|\xi|^2-\frac{1}{4}}\right)-W^1_b(\xi,t)
\right|
\le Ct^b |\xi|^{-b},
\qquad
t\ge1, 
\end{equation*}
on the support of $\chi_H$. 
Since $-2b+n<0$, we have
\begin{equation*}
\begin{split}
\left\|
\chi_H(|\xi|)
\left[
\cos\left(t\sqrt{|\xi|^2-\frac{1}{4}}\right)-W^1_b(\xi,t)
\right]
\right\|_2
& \le Ct^b \||\xi|^{-b}\|_{L^2(|\xi|\ge1)} \\
& \le Ct^b,
\qquad 
t\ge1.
\end{split}
\end{equation*}
On the other hand, we can conclude from Lemma~\ref{lem:H0} that 
\begin{equation}
\label{ineq:p}
\begin{split}
\left\|
\chi_H(|\xi|)D^1_\ell (\xi,t)
\right\|_2 
& \le C\left\|e^{-t|\xi|^2}\right\|_{L^2(|\xi|\ge1)}
+\left\|e^{-t|\xi|^2}\sum_{k=1}^{\ell-1} |\xi|^{2k}
\sum_{j_k=1}^k (t|\xi|^2)^{j_k}
\right\|_{L^2(|\xi|\ge1)} \\
& \le Ct^{\ell-1} \||\xi|^{4(\ell-1)}e^{-t|\xi|^2}\|_{L^2(|\xi|\ge1)} \\
& \le Ct^{\ell-1} \left(e^{-\frac{3}{2}t}\int_{|\xi|\ge1}|\xi|^{8(\ell-1)} e^{-\frac{1}{2}|\xi|^2}\,d\xi\right)^{1/2} \\
& \le Ct^{\ell-1} e^{-\frac{3}{4}t},
\qquad
t\ge1,
\end{split}
\end{equation}
in the case $\ell\ge2$. 
The same estimate also holds for $\ell=1$.
$\Box$
\begin{proposition}
\label{prop:2}
Let $n\in\textbf{N}$, $b\in\textbf{N}$ with $b>n/2$, and $\ell\in\textbf{N}$. 
Then there exists a constant $C>0$ such that 
\begin{equation*}
\left\|m^2_{b,\ell,H}(t)\right\|_2
\le Ct^{b-1} e^{-\frac{t}{2}}
\end{equation*}
for $t\ge1$.
\end{proposition}
\textbf{Proof.} 
First, one considers the case $b\ge2$.
By the Taylor theorem one has 
\begin{equation*}
\begin{split}
\frac{\sin\left(t\sqrt{|\xi|^2-\frac{1}{4}}\right)}{\sqrt{|\xi|^2-\frac{1}{4}}}-W^2_b(\xi,t)
& =2t^{-1}\left[\frac{\partial f}{\partial c}\left(|\xi|,\frac{1}{4},t\right)
-\sum_{k=0}^{b-2} \frac{1}{k!}\left(\frac{1}{4}\right)^k \frac{\partial^{k+1} f}{\partial c^{k+1}}\left(|\xi|,0,t\right)\right] \\
& =\frac{2}{(b-1)!}\left(\frac{1}{4}\right)^{b-1} t^{-1}\frac{\partial^b f}{\partial c^b}\left(|\xi|,\frac{1}{4}\tau,t \right)
\end{split}
\end{equation*}
for some $0\le\tau\le1$. It follows from \eqref{ineq:M2} that 
\begin{equation*}
\left|
\frac{\sin\left(t\sqrt{|\xi|^2-\frac{1}{4}}\right)}{\sqrt{|\xi|^2-\frac{1}{4}}}-W^2_b(\xi,t)
\right|
\le Ct^{b-1} |\xi|^{-b},
\qquad 
t\ge1,
\end{equation*}
on the support of $\chi_H$. 
Thus 
\begin{equation*}
\begin{split}
\left\|
\chi_H(|\xi|)
\left[
\frac{\sin\left(t\sqrt{|\xi|^2-\frac{1}{4}}\right)}{\sqrt{|\xi|^2-\frac{1}{4}}}-W^2_b(\xi,t)
\right]
\right\|_2
& \le Ct^{b-1},
\qquad
t\ge1,
\end{split}
\end{equation*}
in the case $b\ge2$. 
If $b=1$ and $n=1$, then one can obtain 
\begin{equation*}
\begin{split}
\left\|
\chi_H(|\xi|)
\frac{\sin\left(t\sqrt{|\xi|^2-\frac{1}{4}}\right)}{\sqrt{|\xi|^2-\frac{1}{4}}}
\right\|_2 
& \le C\||\xi|^{-1}\|_{L^2(|\xi|\ge1)}
\le C,
\qquad
t>0.
\end{split}
\end{equation*}
On the other hand, from Lemma~\ref{lem:H1}, 
a similar estimate to \eqref{ineq:p} shows that 
\begin{equation*}
\begin{split}
\left\|
\chi_H(|\xi|)D^2_\ell (\xi,t)
\right\|_2 
& \le Ct^{\ell-1} e^{-\frac{3}{4}t},
\qquad
t\ge1,
\end{split}
\end{equation*}
for $\ell\in\textbf{N}$. 
$\Box$
\\

Next we state some estimates in the low frequency region.  
They are the most difficult parts to be proved in this paper. 
\begin{proposition}
\label{prop:3}
Let $n\in\textbf{N}$, $b\in\textbf{N}$ and $\ell\in\textbf{N}$. 
Then there exists a constant $C>0$ such that 
\begin{equation*}
\left\|m^1_{b,\ell,L}(t)\right\|_2
\le Ct^{2(b-1)}e^{-\frac{t}{2}}
+Ct^{-\frac{n}{4}-\ell} 
\end{equation*}
for $t\ge1$.
\end{proposition}
\textbf{Proof.} 
If $|\xi|\le1/2$, then 
\begin{equation*}
\begin{split}
e^{-\frac{t}{2}}\cosh\left(t\sqrt{\frac{1}{4}-|\xi|^2}\right) 
& =\frac{1}{2}e^{-t\left(\frac{1}{2}-\sqrt{\frac{1}{4}-|\xi|^2}\right)}
+\frac{1}{2}e^{-t\left(\frac{1}{2}+\sqrt{\frac{1}{4}-|\xi|^2}\right)}.
\end{split}
\end{equation*}
By definition of $g$, we have $e^{-t\left(\frac{1}{2}-\sqrt{\frac{1}{4}-|\xi|^2}\right)} 
=g(|\xi|,1,t)$, since 
\begin{equation*}
\begin{split}
\frac{1}{2}-\sqrt{\frac{1}{4}-|\xi|^2} 
=\frac{|\xi|^2}{\frac{1}{2}+\sqrt{\frac{1}{4}-|\xi|^2}} 
=\frac{2|\xi|^2}{1+\sqrt{1-4|\xi|^2}}.
\end{split}
\end{equation*}
Hence, $m_{b,\ell,L}^1$ can be represented as  
\begin{equation*}
\begin{split}
m^1_{b,\ell,L}(\xi,t) 
& =\chi_L(|\xi|)\left[
\biggr(
\frac{1}{2} g(|\xi|,1,t)
-D^1_\ell(\xi,t)
\biggr)
+\frac{1}{2}e^{-t\left(\frac{1}{2}+\sqrt{\frac{1}{4}-|\xi|^2}\right)}
-e^{-\frac{t}{2}}
W_b^1(\xi,t)
\right].
\end{split}
\end{equation*}
By the Taylor theorem we have 
\begin{equation*}
\frac{1}{2} g(|\xi|,1)
-D^1_\ell(\xi)
=\frac{1}{2\ell!}\frac{\partial^\ell g}{\partial a^\ell}(|\xi|,\tau,t)
\end{equation*}
for some $0\le\tau\le1$.
Therefore, it follows from \eqref{ineq:H0} and Lemma~\ref{lem:N2} that 
\begin{equation*}
\label{ineq:p2}
\begin{split}
\left\|
\chi_L(|\xi|)\biggr(
\frac{1}{2} g(|\xi|,1,t)
-D^1_\ell(\xi,t)
\biggr)
\right\|_2 
& \le C\left\|
|\xi|^{2\ell}e^{-t|\xi|^2} \sum_{j=1}^\ell (t|\xi|^2)^j 
\right\|_{L^2(|\xi|\le1)} \\
& \le C\sum_{j=1}^\ell 
t^j \left\|
|\xi|^{2(j+\ell)}e^{-t|\xi|^2}
\right\|_{L^2(|\xi|\le1)} \\
& \le C\sum_{j=1}^\ell 
t^j (1+t)^{-\frac{n}{4}-(j+\ell)} \\
& \le C(1+t)^{-\frac{n}{4}-\ell},
\qquad
t>0,
\end{split}
\end{equation*}
for $\ell\in\textbf{N}$.
Furthermore, we find that 
\begin{equation*}
\begin{split}
\left\|
\chi_L(|\xi|)e^{-t\left(\frac{1}{2}+\sqrt{\frac{1}{4}-|\xi|^2}\right)}
\right\|_2
\le Ce^{-t(\frac{1}{2}+\frac{\sqrt{5}}{6})},
\qquad
t>0.
\end{split}
\end{equation*}
On the other hand, \eqref{eq:SING} and \eqref{ineq:M3} imply 
\begin{equation*}
\begin{split}
\left\|
\chi_L(|\xi|)
e^{-\frac{t}{2}}
W_b^1(\xi,t)\right\|_2
 \le Ce^{-\frac{t}{2}}\sum_{k=0}^{b-1} t^{2k} 
 \le Ct^{2(b-1)}e^{-\frac{t}{2}}, 
\qquad 
t\ge1.
\end{split}
\end{equation*}
Combining above three inequalities, one obtains the desired estimates. 
$\Box$
\begin{proposition}
Let $n\in\textbf{N}$, $b\in\textbf{N}$ and $\ell\in\textbf{N}$. 
Then there exists a constant $C>0$ such that 
\begin{equation*}
\left\|m^2_{b,\ell,L}(t)\right\|_2
\le 
\begin{cases}
Ct^{2(b-1)-1}e^{-\frac{t}{2}}+Ct^{-\frac{n}{4}-\ell}, & b\ge2, \\
Ce^{-\frac{t}{2}}+Ct^{-\frac{n}{4}-\ell}, & b=1,
\end{cases}
\end{equation*}
for $t\ge1$.
\end{proposition}
\textbf{Proof.} 
Similarly to the proof of Proposition~\ref{prop:3}, we estimate  
\begin{equation*}
\begin{split}
m^2_{b,\ell,L}(\xi) 
& =\chi_L(|\xi|)\left[
\biggr(
h(|\xi|,1,t)
-D^2_\ell(\xi,t)
\biggr)
-\frac{e^{-t\left(\frac{1}{2}+\sqrt{\frac{1}{4}-|\xi|^2}\right)}}{\sqrt{1-4|\xi|^2}}
-e^{-\frac{t}{2}}
W_b^2(\xi,t)
\right].
\end{split}
\end{equation*}
By the Taylor theorem we have 
\begin{equation*}
h(|\xi|,1,t)
-D^2_\ell(\xi,t)
=\frac{1}{\ell!}\frac{\partial^\ell h}{\partial a^\ell}(|\xi|,\tau,t)
\end{equation*}
for some $0\le\tau\le1$.
We use \eqref{ineq:H1} and Lemma~\ref{lem:N2} to have  
\begin{equation*}
\begin{split}
\left\|
\chi_L(|\xi|)\biggr(
h(|\xi|,1,t)
-D^2_\ell(\xi,t)
\biggr)
\right\|_2 
& \le C\left\|
|\xi|^{2\ell}e^{-t|\xi|^2} \sum_{j=0}^\ell (t|\xi|^2)^j 
\right\|_{L^2(|\xi|\le1)} \\
& \le C(1+t)^{-\frac{n}{4}-\ell},
\qquad
t>0, 
\end{split}
\end{equation*}
for $\ell\in\textbf{N}$.
Next, it can be easily seen that 
\begin{equation*}
\begin{split}
\left\|
\chi_L(|\xi|)
\frac{e^{-t\left(\frac{1}{2}+\sqrt{\frac{1}{4}-|\xi|^2}\right)}}{\sqrt{1-4|\xi|^2}}
\right\|_2
\le Ce^{-t(\frac{1}{2}+\frac{\sqrt{5}}{6})},
\qquad
t>0.
\end{split}
\end{equation*}
Furthermore, by \eqref{eq:SING}  and \eqref{ineq:M3}, one has 
\begin{equation*}
\begin{split}
\left\|
\chi_L(|\xi|)
e^{-\frac{t}{2}}
W_b^2(\xi,t)\right\|_2 
& \le Ce^{-\frac{t}{2}}t^{-1}
\sum_{k=0}^{b-2}
\left\|
\frac{\partial^{k+1} f}{\partial c^{k+1}}\left(|\xi|,0,t\right)
\right\|_{L^2(|\xi|\le1)} \\
& \le Ce^{-\frac{t}{2}}t^{-1}\sum_{k=0}^{b-2} t^{2(k+1)} \\
& \le Ct^{2(b-1)-1}e^{-\frac{t}{2}},
\qquad 
t\ge1, 
\end{split}
\end{equation*}
in the case $b\ge2$. 
If $b=1$, then $W_b^2(\xi,t)=0$, and the proof is now complete.
$\Box$
\\

Finally, we deal with the decay estimates in the middle frequency region. 
It is much easier to obtain the following estimates.   
\begin{proposition}
\label{prop:5}
Let $n\in\textbf{N}$, $b\in\textbf{N}$ and $\ell\in\textbf{N}$. 
Then there exists a constant $C>0$ such that 
\begin{equation*}
\label{ineq:mid}
\left\|m^1_{b,\ell,M}(t)\right\|_2
\le Ct^{b-1}e^{-\frac{t}{2}}
+Ct^{\ell-1}e^{-\frac{t}{9}}
\end{equation*}
for $t\ge1$.
\end{proposition}
\textbf{Proof.} 
According to \eqref{ineq:M3} and \eqref{eq:H0}, 
we obtain 
\begin{equation*}
\begin{split}
& \left\|m^1_{b,\ell,M}(t)\right\|_2 \\
& \le e^{-\frac{t}{2}}
+Ce^{-\frac{t}{2}}
\sum_{k=0}^{b-1}\left\|
\frac{\partial^k f}{\partial c^k}\left(|\xi|,0,t\right)
\right\|_{L^2(1/3\le|\xi|\le1)}
+C\sum_{k=0}^{\ell-1} \left\|
\frac{\partial^k g}{\partial a^k}(|\xi|,0,t)
\right\|_{L^2(1/3\le|\xi|\le1)} \\
& \le e^{-\frac{t}{2}}
+Ce^{-\frac{t}{2}}
\sum_{k=0}^{b-1}
t^k\||\xi|^{-k}\|_{L^2(1/3\le|\xi|\le1)}
+C\sum_{k=0}^{\ell-1}
t^k \|(1+|\xi|^2)^{2k}e^{-t|\xi|^2}\|_{L^2(1/3\le|\xi|\le1)} \\
& \le Ct^{b-1}e^{-\frac{t}{2}}
+Ct^{\ell-1}e^{-\frac{t}{9}}
\end{split}
\end{equation*}
for $t\ge1$. 
$\Box$
\begin{proposition}
Let $n\in\textbf{N}$, $b\in\textbf{N}$ and $\ell\in\textbf{N}$. 
Then there exists a constant $C>0$ such that 
\begin{equation}
\label{ineq:mid2}
\left\|m^2_{b,\ell,M}(t)\right\|_2
\le Ct^{b-1}e^{-\frac{t}{2}}
+Ct^{\ell-1}e^{-\frac{t}{9}}
\end{equation}
for $t\ge1$.
\end{proposition}
\textbf{Proof.} 
Similarly to the proof of Proposition~\ref{prop:5}, 
it follows from \eqref{ineq:M3} and \eqref{eq:H1} that 
\begin{equation*}
\begin{split}
& \left\|m^2_{b,\ell,M}(t)\right\|_2 \\
& \le te^{-\frac{t}{2}}
+Ce^{-\frac{t}{2}}
\sum_{k=0}^{b-2} \left\|
\frac{\partial^{k+1} f}{\partial c^{k+1}}(|\xi|,0,t)
\right\|_{L^2(1/3\le|\xi|\le1)} 
+C\sum_{k=0}^{\ell-1}
\left\|
\frac{\partial^k h}{\partial a^k}(|\xi|,0,t)
\right\|_{L^2(1/3\le|\xi|\le1)} \\
& \le te^{-\frac{t}{2}}
+Ce^{-\frac{t}{2}}
\sum_{k=1}^{b-1} 
t^k \||\xi|^{-k}\|_{L^2(1/3\le|\xi|\le1)} 
+C\sum_{k=0}^{\ell-1}
t^k \|(1+|\xi|^2)^{2k}e^{-t|\xi|^2}\|_{L^2(1/3\le|\xi|\le1)} \\
& \le Ct^{b-1}e^{-\frac{t}{2}}
+Ct^{\ell-1}e^{-\frac{t}{9}},
\qquad 
t\ge1,
\end{split}
\end{equation*}
in the case $b\ge2$. When $b=1$, it is much simpler to prove \eqref{ineq:mid2}. 
$\Box$

\section{Appendix}
In this appendix, we check the equivalence of the representations between $D_\ell^i$ ($i=1$, $2$) and that derived by the Takeda expansion. 

We first check that our results include those of \cite{HO}, \cite{MN} and \cite{Ni}. 
For this, we list several derivatives of $f$, $g$ and $h$ in order to see the asymptotic profile of the solution to \eqref{eq:LDW}. 
For simplicity, we write $r=|\xi|$ throughout this section. 

The derivatives of $f$ are
\begin{equation*}
\frac{\partial f}{\partial c}(r,0,t)
=\frac{t\sin(tr)}{2r}, 
\qquad
\frac{\partial^2 f}{\partial^2 c}(r,0,t)
=\frac{t\sin(tr)-t^2 r\cos(tr)}{4r^3},
\end{equation*}
\begin{equation*}
\frac{\partial^3 f}{\partial^3 c}(r,0,t)
=\frac{3t\sin(tr)-3t^2 r \cos(tr)-t^3 r^2 \sin(tr)}{8r^5}.
\end{equation*}
This shows \eqref{as} and \eqref{as2} from Theorem~\ref{thm:1} with $b=2$, $\ell=1$ and $b=1$, $\ell=1$, respectively, which includes the previous results of \cite{HO}, \cite{MN}, and \cite{Ni}. 
In the previous studies, the asymptotic profile of wave part was not found in the case $n=4$. 
However, by substituting $b=3$ into $W_b^i(r,t)$ ($i=1$, $2$), we obtain  
\begin{equation*}
\begin{split}
& e^{-\frac{t}{2}}
\left(
\cos(tr)
+\frac{t\sin(tr)}{8r}
+\frac{t\sin(tr)-t^2 r\cos(tr)}{128 r^3}
\right)\widehat{u_0} \\
& \qquad\qquad
+e^{-\frac{t}{2}}
\left(
\frac{\sin(tr)}{r}
+\frac{\sin(tr)-t r\cos(tr)}{8r^3}
\right)\left(\frac{1}{2}\widehat{u_0}+\widehat{u_1}\right),
\end{split}
\end{equation*}
which is completely new. 

Next, the derivatives of $g$ are  
\begin{equation*}
\frac{\partial g}{\partial a}(r,0,t)
=r^2 (-tr^2)e^{-tr^2}, 
\qquad
\frac{\partial^2 g}{\partial^2 a}(r,0,t)
=4r^4 (-tr^2)e^{-tr^2}
+r^4 (-tr^2)^2 e^{-tr^2}.
\end{equation*}
and the derivatives of $h$ are 
\begin{equation*}
\frac{\partial h}{\partial a}(r,0,t)
=2r^2 e^{-tr^2}
+r^2 (-tr^2)e^{-tr^2}, 
\end{equation*}
\begin{equation*}
\frac{\partial^2 h}{\partial a^2}(r,0,t)
=12r^4 e^{-tr^2}
+8r^4 (-tr^2) e^{-tr^2}
+r^4 (-tr^2)^2 e^{-tr^2}.
\end{equation*}

Now, let us check the equivalence between the asymptotic profiles obtained in Theorem~\ref{thm:1} and those of the Takeda expansion. 

We first calculate coefficients defined by \eqref{coe} of the Takeda expansion:
\begin{equation*}
\begin{split}
& \alpha_{0,0}=1, \\
& \alpha_{0,1}=0, 
\quad 
\alpha_{1,0}=1, \\
& \alpha_{0,2}=0,
\quad
\alpha_{1,1}=2, 
\quad 
\alpha_{2,0}=1/2, \\
& \alpha_{0,3}=0,
\quad
\alpha_{1,2}=5, 
\quad 
\alpha_{2,1}=2, 
\quad 
\alpha_{3,0}=1/6, \\
& \beta_0=1,
\qquad 
\beta_1=2,
\qquad 
\beta_2=6,
\qquad 
\beta_3=20.
\end{split}
\end{equation*}
On the coefficients $\beta_\ell$, in general, we see that 
\begin{equation*}
\beta_\ell
=\frac{2^\ell (2\ell-1)!!}{\ell!},
\qquad
\ell\in\textbf{N},
\end{equation*}
and the quantity $2^\ell (2\ell-1)!!$ has been already seen in \eqref{eq:SING} 
(see also \eqref{eq:const}, \eqref{eq:2con} and \eqref{eq:cco}).
Thus the first order expansion is
\begin{equation*}
\left\{
e^{-tr^2}
+r^2(-tr^2)e^{-tr^2}
\right\}\widehat{u_0}
+\left\{\left(1+2 r^2\right)e^{-tr^2}
+r^2 (-tr^2) e^{-tr^2}
\right\}
\left(\frac{1}{2}\widehat{u_0}+\widehat{u_1}\right).
\end{equation*}
and the second order expansion is 
\begin{equation*}
\begin{split}
& \left\{
e^{-tr^2}
+(r^2+2r^4)(-tr^2)e^{-tr^2}
+\frac{1}{2}r^4 (-tr^2)^2 e^{-tr^2}
\right\}\widehat{u_0} \\
& \qquad
+\left\{\left(1+2 r^2+6r^4\right)e^{-tr^2}
+(r^2+4r^4) (-tr^2) e^{-tr^2}
+\frac{1}{2}r^4 (-tr^2)^2 e^{-tr^2}
\right\}
\left(\frac{1}{2}\widehat{u_0}+\widehat{u_1}\right).
\end{split}
\end{equation*}
Compared with the author's expansion (substituting $\ell=1$ and $\ell=2$ into $D_\ell^i$, $i=1$, $2$), 
we can find that both representations are the same. 
The following two theorems assure that they are completely the same expressions not only in the lower order case but also for the higher order expansions.
\begin{theorem}
\label{thm:E1}
For each $m\in\textbf{N}$, 
\begin{equation}
\label{eq:EE}
\sum_{k=0}^m \frac{1}{k!}
\frac{\partial^k g}{\partial a^k}(r,0)
=\sum_{j=0}^m \sum_{k=0}^{m-j} \left.\frac{1}{j! k!}\frac{d^k}{ds^k}\phi_1(s)^j\right|_{s=0}(-tr^2)^j r^{2(j+k)}e^{-tr^2},
\end{equation}
where 
\begin{equation*}
g(r,a)=\exp\left(-\frac{2tr^2}{1+\sqrt{1-4ar^2}}\right),
\qquad
\phi_1(s)=\left(\frac{1}{1/2+\sqrt{1/4-s}}\right)^2.
\end{equation*}
\end{theorem}
\textbf{Proof.} 
Since  
\begin{equation*}
\begin{split}
\frac{1}{\frac{1}{2}+\sqrt{\frac{1}{4}-ar^2}}-1 
& =\frac{\frac{1}{2}-\sqrt{\frac{1}{4}-ar^2}}{\frac{1}{2}+\sqrt{\frac{1}{4}-ar^2}} 
=ar^2 \phi_1(ar^2),
\end{split}
\end{equation*}
we see that 
\begin{equation*}
\begin{split}
e^{tr^2}g(r,a)
& =\exp\left(tr^2-\frac{2tr^2}{1+\sqrt{1-4ar^2}}\right) 
=\exp\left\{-tr^2\left(\frac{1}{\frac{1}{2}+\sqrt{\frac{1}{4}-ar^2}}-1\right)\right\} \\
& =\exp\biggr(-tr^2\cdot ar^2\phi_1(ar^2)\biggr).
\end{split}
\end{equation*}
Moreover, it follows that 
\begin{equation*}
\frac{\partial^k}{\partial a^k}\exp\biggr(-tr^2\cdot ar^2\phi_1(ar^2)\biggr)
=r^{2k} \left.
\frac{\partial^k}{\partial s^k}\exp\biggr(-tr^2\cdot s\phi_1(s)\biggr)
\right|_{s=ar^2}.
\end{equation*}
On the other hand, we have 
\begin{equation*}
\begin{split}
\sum_{j=0}^m \sum_{k=0}^{m-j} \left.\frac{1}{j! k!}\frac{d^k}{ds^k}\phi_1(s)^j\right|_{s=0}(-tr^2)^j r^{2(j+k)}
& =\sum_{k=0}^m \sum_{j=0}^k 
\left.\frac{1}{j!(k-j)!}\frac{d^{k-j}}{ds^{k-j}}\phi_1(s)^j\right|_{s=0}
(-tr^2)^j r^{2k}.
\end{split}
\end{equation*}
Thus \eqref{eq:EE} is equivalent to 
\begin{equation}
\label{eq:EE2}
\sum_{k=0}^m \frac{1}{k!}
\left.
\frac{\partial^k}{\partial s^k}\exp\biggr(-tr^2\cdot s\phi_1(s)\biggr)
\right|_{s=0}
=\sum_{k=0}^m \sum_{j=0}^k 
\left.\frac{1}{j!(k-j)!}\frac{d^{k-j}}{ds^{k-j}}\phi_1(s)^j\right|_{s=0}
(-tr^2)^j.
\end{equation}
Putting $\phi(s):=-tr^2 \phi_1(s)$, \eqref{eq:EE2} becomes 
\begin{equation*}
\sum_{k=0}^m \frac{1}{k!}
\left.
\frac{\partial^k}{\partial s^k}\exp\biggr(s\phi(s)\biggr)
\right|_{s=0}
=\sum_{k=0}^m \sum_{j=0}^k 
\left.\frac{1}{j!(k-j)!}\frac{d^{k-j}}{ds^{k-j}}\phi(s)^j\right|_{s=0}.
\end{equation*}
So it suffices to show that 
\begin{equation}
\label{eq:ExP}
\left.
\frac{\partial^k}{\partial s^k}\exp\biggr(s\phi(s)\biggr)
\right|_{s=0}
=\sum_{j=0}^k 
\left.\frac{k!}{j!(k-j)!}\frac{d^{k-j}}{ds^{k-j}}\phi(s)^j\right|_{s=0}.
\end{equation}
However, by the Taylor theorem one has the equality 
\begin{equation*}
\begin{split}
\left.\frac{d^k}{ds^k}\exp\biggr(s\phi(s)\biggr)\right|_{s=0} 
& =\left.\left(\frac{d^k}{ds^k}\sum_{j=0}^\infty \frac{1}{j!}\biggr(s\phi(s)\biggr)^j\right)\right|_{s=0} \\
& =\left.\sum_{j=0}^\infty \frac{1}{j!}\frac{d^k}{ds^k}\biggr(s\phi(s)\biggr)^j\right|_{s=0} 
=\left.\sum_{j=0}^k \frac{1}{j!}\frac{d^k}{ds^k}\biggr(s\phi(s)\biggr)^j\right|_{s=0} \\
& =\sum_{j=0}^k \frac{1}{j!}
\sum_{i=0}^k \frac{k!}{i!(k-i)!}
\left.\left(\frac{d^i}{ds^i}s^j\right)\right|_{s=0}
\left.\left(\frac{d^{k-i}}{ds^{k-i}}\phi(s)^j\right)\right|_{s=0} \\
& =\sum_{j=0}^k \frac{1}{j!}
\frac{k!}{j!(k-j)!}
j!
\left.\left(\frac{d^{k-j}}{ds^{k-j}}\phi(s)^j\right)\right|_{s=0} \\
& =\sum_{j=0}^k
\frac{k!}{j!(k-j)!}
\left.\left(\frac{d^{k-j}}{ds^{k-j}}\phi(s)^j\right)\right|_{s=0}
\end{split}
\end{equation*}
for any $k\in\textbf{N}$, which implies the validity of \eqref{eq:ExP}. 
$\Box$
\begin{theorem}
For each $m\in\textbf{N}$,
\begin{equation*}
\label{eq:SEE}
\begin{split}
\sum_{k=0}^m \frac{1}{k!}\frac{\partial^k h}{\partial a^k}(r,0) 
=\sum_{j=0}^m \sum_{k=0}^{m-j}\sum_{\ell=0}^{m-j-k}
\frac{1}{j!k!\ell!}
\left.\frac{d^k}{ds^k}\phi_1(s)^j\right|_{s=0}
\left.\frac{d^k}{ds^k}\psi(s)\right|_{s=0}
(-tr^2)^j r^{2(j+k+\ell)}e^{-tr^2},
\end{split}
\end{equation*}
where 
\begin{equation*}
\phi_1(s)=\left(\frac{1}{1/2+\sqrt{1/4-s}}\right)^2,
\qquad
\psi(s)=\frac{1}{\sqrt{1-4r}},
\end{equation*}
\begin{equation*}
h(r,a)=\frac{1}{\sqrt{1-4ar^2}}\exp\left(-\frac{2tr^2}{1+\sqrt{1-4ar^2}}\right).
\end{equation*}
\end{theorem}
\textbf{Proof.} 
Similarly to the proof of Theorem~\ref{thm:E1}, 
it suffices to show that 
\begin{equation}
\label{eq:ExP2}
\begin{split}
\left.
\frac{d^k}{ds^k}\left(\psi(s)\exp\biggr(s\phi(s)\biggr)\right)
\right|_{s=0} 
=\sum_{j=0}^k \sum_{\ell=0}^{k-j}
\frac{k!}{j!(k-j-\ell)!\ell!}
\left.\frac{d^{k-j-\ell}}{ds^{k-j-\ell}}\phi(s)^j\right|_{s=0}
\left.\frac{d^\ell}{ds^\ell}\psi(s)\right|_{s=0},
\end{split}
\end{equation}
since $h(r,a)=\psi(ar^2)g(r,a)$.
In fact, it follows that 
\begin{equation*}
\begin{split}
& \left.
\frac{d^k}{ds^k}\left(\psi(s)\exp\biggr(s\phi(s)\biggr)\right)
\right|_{s=0} 
=\sum_{j=0}^k \frac{1}{j!} 
\left.
\frac{d^k}{ds^k}\left(\biggr(s\phi(s)\biggr)^j \psi(s)\right)\right|_{s=0} \\
& =\sum_{j=0}^k \frac{1}{j!} 
\sum_{\ell=0}^{k} \frac{k!}{(k-\ell)!\ell!}
\left.\frac{d^{k-\ell}}{ds^{k-\ell}}\biggr(s\phi(s)\biggr)^j\right|_{s=0}
\left.\frac{d^\ell}{ds^\ell}\psi(s)\right|_{s=0} \\
& =\sum_{j=0}^k \frac{1}{j!} 
\sum_{\ell=0}^{k-j} \frac{k!}{(k-\ell)!\ell!}
\left.\frac{d^{k-\ell}}{ds^{k-\ell}}\biggr(s\phi(s)\biggr)^j\right|_{s=0}
\left.\frac{d^\ell}{ds^\ell}\psi(s)\right|_{s=0} \\
& =\sum_{j=0}^k \frac{1}{j!} 
\sum_{\ell=0}^{k-j} \frac{k!}{(k-\ell)!\ell!}
\frac{(k-\ell)!}{(k-\ell-j)!}\left.\frac{d^{k-\ell-j}}{ds^{k-\ell-j}}\phi(s)^j\right|_{s=0}
\left.\frac{d^\ell}{ds^\ell}\psi(s)\right|_{s=0} \\
& =\sum_{j=0}^k \sum_{\ell=0}^{k-j}
\frac{k!}{j!(k-j-\ell)!\ell!}
\left.\frac{d^{k-j-\ell}}{ds^{k-j-\ell}}\phi(s)^j\right|_{s=0}
\left.\frac{d^\ell}{ds^\ell}\psi(s)\right|_{s=0}
\end{split}
\end{equation*}
for any $k\in\textbf{N}$. 
$\Box$

\begin{remark}
One can obtain \eqref{eq:ExP2} also from \eqref{eq:ExP} {\rm :}
\begin{equation*}
\begin{split}
& \left.
\frac{d^k}{ds^k}\left(\psi(s)\exp\biggr(s\phi(s)\biggr)\right)
\right|_{s=0} \\
& =\sum_{\ell=0}^k \frac{k!}{(k-\ell)!\ell!}
\left.\frac{d^{k-\ell}}{ds^{k-\ell}}\exp\biggr(s\phi(s)\biggr)\right|_{s=0}
\left.\frac{d^\ell}{ds^\ell}\psi(s)\right|_{s=0} \\
& =\sum_{\ell=0}^k \frac{k!}{(k-\ell)!\ell!}
\sum_{j=0}^{k-\ell}\frac{(k-\ell)!}{j!(k-\ell-j)!}
\left.\frac{d^{k-\ell-j}}{ds^{k-\ell-j}}\phi(s)^j\right|_{s=0}
\left.\frac{d^\ell}{ds^\ell}\psi(s)\right|_{s=0} \\
& =\sum_{\ell=0}^k \sum_{j=0}^{k-\ell}
\frac{k!}{j!(k-\ell-j)!\ell!}
\left.\frac{d^{k-\ell-j}}{ds^{k-\ell-j}}\phi(s)^j\right|_{s=0}
\left.\frac{d^\ell}{ds^\ell}\psi(s)\right|_{s=0} \\
& =\sum_{j=0}^k \sum_{\ell=0}^{k-j}
\frac{k!}{j!(k-j-\ell)!\ell!}
\left.\frac{d^{k-j-\ell}}{ds^{k-j-\ell}}\phi(s)^j\right|_{s=0}
\left.\frac{d^\ell}{ds^\ell}\psi(s)\right|_{s=0}.
\end{split}
\end{equation*}
\end{remark}
%

{\bf Acknowledgment.} 
The author would like to thank Professor Ryo Ikehata for valuable discussion and kind encouragement. 
The author is also grateful to Professor Kenji Nishihara for useful comment. 
The author is also grateful to Mr. Kenji Kurogi and Mr. Go Nakamura for helpful suggestions.

\bibliographystyle{amsplain}

\begin{thebibliography}{99}

\bibitem{CH}
R. Courant and D. Hilbert, 
Methods of Mathematical Physics, Vol. II, Wiley, New York, 1962.

\bibitem{EZ}
M. Escobedo and E. Zuazua, 
\textit{Large time behavior for convection-diffusion equation in $R^n$}, 
J. Funct. Anal. \textbf{100} (1991), 119--161.

\bibitem{F}
H. Fujita, 
\textit{On the blowing up of solutions of the Cauchy problem for $u_t-\Delta u=u^{1+\alpha}$}, 
J. Fac. Sci. Univ. Tokyo \textbf{13} (1966), 109--124.

\bibitem{HKN}
N. Hayashi, E. Kaikina and P.I. Naumkin, 
\textit{Damped wave equation with super critical nonlinearities}, 
Differential Integral Equations \textbf{17} (2004), 637--652

\bibitem{HKN2}
N. Hayashi, E. Kaikina and P.I. Naumkin, 
\textit{Damped wave equation with a critical nonlinearity}, 
Trans. Amer. Math. Soc. \textbf{358} (2006), 1165--1185.

\bibitem{HL}
L. Hsiao and T.-P. Liu,  
\textit{Convergence to nonlinear diffusion waves for solutions of a system of hyperbolic conservations with damping}, 
Comm. Math. Phys. \textbf{143} (1992), 599--605.

\bibitem{HL2}
L. Hsiao and T.-P. Liu, 
\textit{Nonlinear diffusive phenomena of nonlinear hyperbolic systems}, 
Chin. Ann. Math. Ser. B \textbf{14} (1993), 465--480

\bibitem{HO}
T. Hosono and T. Ogawa, 
\textit{Large time behavior and $L^p$-$L^q$ estimate of $2$-dimensional nonlinear damped wave equations}, 
J. Differential Equations \textbf{203} (2004), 82--118.

\bibitem{IIK}
K. Ishige, M. Ishiwata and T. Kawakami, 
\textit{The decay of the solutions for the heat equation with a potential}, 
Indiana Univ. Math. J. \textbf{58} (2009), 2673--2707. 

\bibitem{IKK}
K. Ishige, T. Kawakami and K. Kobayashi, 
\textit{Asymptotics for a nonlinear integral equation with a generalized heat kernel}, 
J. Evol. Equ. \textbf{14} (2014), 749--777.

\bibitem{IKM}
K. Ishige, T. Kawakami and  H. Michihisa, 
\textit{Asymptotic expansions of solutions of fractional diffusion equations}, 
SIAM J. Math. Anal. \textbf{49} (2017), 2167--2190.

\bibitem{K}
G. Karch, 
\textit{Selfsimilar profiles in large time asymptotics of solutions to damped wave equations}, 
Studia Math. \textbf{143} (2000), 175--197.

\bibitem{KNO}
S. Kawashima, M. Nakao and K. Ono,
\textit{On the decay property of solutions to the Cauchy problem of the semilinear wave equation with a dissipative term}, 
J. Math. Soc. Japan \textbf{47} (1995), 617--653.

\bibitem{IMN}
R. Ikehata, Y. Miyaoka and T. Nakatake, 
\textit{Decay estimates of solutions for dissipative wave equations in $R^N$ with lower power nonlinearities}, 
J. Math. Soc. Japan \textbf{56} (2004), 365--373. 

\bibitem{IT}
R. Ikehata and K. Tanizawa, 
\textit{Global existence of solutions for semilinear damped wave equations in $R^N$ with noncompactly supported initial data}, 
Nonlinear Anal. \textbf{61} (2005), 1189--1208.

\bibitem{MN}
P. Marcati and K. Nishihara, 
\textit{The $L^p$-$L^q$ estimates of solutions to one-dimensional damped wave equations and their application to compressible flow through porous media}, 
J. Differential Equations \textbf{191} (2003), 445--469.

\bibitem{M}
A. Matsumura, 
\textit{On the asymptotic behavior of solutions of semilinear wave equations}, 
Publ. Res. Inst. Sci. Kyoto Univ. \textbf{12} (1976), 169--189.

\bibitem{Na}
T. Narazaki, 
\textit{$L^p$-$L^q$ estimates for damped wave equations and their applications to semi-linear problem}, 
J. Math. Soc. Japan \textbf{56} (2004), 585--626.

\bibitem{Ni}
K. Nishihara, 
\textit{$L^p$-$L^q$ estimates of solutions to the damped wave equation in $3$-dimensional space and their application}, 
Math. Z. \textbf{244} (2003), 631--649.

\bibitem{SW}
S. Sakata and Y. Wakasugi,
\textit{Movement of time-delayed hot spots in Euclidean space}, 
Math. Z. \textbf{285} (2017), 1007--1040.

\bibitem{T}
H. Takeda, 
\textit{Higher-order expansion of solutions for a damped wave equation}, 
Asymptot. Anal. \textbf{94} (2015), 1--31.

\bibitem{TY}
H.Takeda and S.Yoshikawa, 
\textit{Asymptotic profiles of solutions for the isothermal Falk-Konopka system of shape memory alloys with weak damping}, 
Asymptotic Anal. \textbf{81} (2013), 331--372. 

\bibitem{ToYo}
G. Todorova and B. Yordanov, 
\textit{Critical exponent for a nonlinear wave equation with damping}, 
J. Differential Equations \textbf{174} (2001), 464--489.

\bibitem{YM}
H. Yang and A. Milani, 
\textit{On the diffusion phenomenon of quasilinear hyperbolic waves}, 
Bull. Sci. Math. \textbf{124} (2000), 415--433.

\end{thebibliography}

\end{document}